\documentclass[final,3p,times]{elsarticle}

\usepackage{graphicx,amsmath,amssymb,txfonts}
\usepackage{algorithm}
\usepackage{algorithmic}
\usepackage{subfigure,color}
\usepackage[colorlinks]{hyperref}
\usepackage{hypernat}
\usepackage{multirow}
\numberwithin{equation}{section}
\graphicspath{{Figs/}}
\usepackage{booktabs}
\usepackage{listings}

\biboptions{sort&compress}

\makeatletter

\newcommand{\Rmnum}[1]{\expandafter\@slowromancap\romannumeral #1@}
\makeatother

\newtheorem{theorem}{\textbf{Theorem}}

\newtheorem{definition}{\textbf{Definition}}
\newtheorem{lemma}{\textbf{Lemma}}

\newtheorem{remark}{\textbf{Remark}}
\newtheorem{example}{\textbf{Example}}

\begin{document}

\begin{frontmatter}



\title{A spatial sixth-order CCD-TVD method for solving multidimensional coupled Burgers' equation}


\author[CSU]{Kejia Pan}
\author[CSU]{Xiaoxin Wu}
\author[XTU]{Xiaoqiang Yue\corref{cor}}
\author[CSU]{Runxin Ni}

\cortext[cor]{Corresponding author.\\E-mail address: yuexq@xtu.edu.cn(X.Q. Yue).}
\address[CSU]{School of Mathematics and Statistics, Central South University, Changsha, Hunan 410083, China}
\address[XTU]{School of Mathematics and Computational Science, Hunan Key Laboratory for Computation
and Simulation in Science and Engineering, Xiangtan University, Hunan 411105, China}

\begin{abstract}
In this paper, an efficient and high-order accuracy finite difference method is proposed for solving multidimensional nonlinear Burgers' equation.
The third-order three stage Runge-Kutta total variation diminishing (TVD) scheme is employed for the time discretization, and the three-point combined compact difference (CCD) scheme is used for spatial discretization. Our method is third-order accurate in time and sixth-order accurate in space.
The CCD-TVD method treats the nonlinear term explicitly thus it is very efficient and easy to implement. In addition, we prove the unique solvability
of the CCD system under non-periodic boundary conditions. Numerical experiments including both two-dimensional and three-dimensional
problems  have been conducted to demonstrate the high efficiency and accuracy of the proposed method.
\end{abstract}


\end{frontmatter}

\section{Introduction}
As an important partial differential equation in fluid mechanics, Burgers'  equation is renowed for its variety of application in traffic flow~\cite{Arpad}, gas dynamics~\cite{Yang}, shock waves~\cite{Khesin}, wave propagation in acoustics~\cite{Davidson}, shallow water waves, etc~\cite{Su,Zabusky,Zhao,Esipov,Logan,Shandarin}.
The analytical solution of nonlinear  Burgers' equation  can be obatined with a special set of  initial and boundary conditions~\cite{Bateman}.
The transformation from the Burgers' equation to a linear heat equation was found by Hopf~\cite{Hopf} and Cole~\cite{Cole} independently, which prompted the study of analytical solutions under the homogeneous boundary conditions.
With the initial vorticity being zero, the exact solutions of the 2D/3D coupled Burgers' equations were obtained using 2D/3D versions of Hopf-Cole transformations~\cite{Fletcher,GaoZou}.  However, as a result of its nonlinear nature,  analytical solutions under the general initial boundary conditions can not be obtained, which makes numerical computations of essential importance.

In the existing literature, it has attracted great efforts to solve the 1D Burgers' equation and its variational form with different numerical methods.
For instance, the finite element methods were developed to solve the nonlinear  Burgers' equation~\cite{Kutluay2013,Finn,Caldwell}.
The mixed finite difference and Galerkin methods were developed by Dehghan et al.~\cite{Dehghan2014} for solving the (coupled) Burgers' equations.
A mixed finite difference and boundary element approach was utilized by Bahadir~\cite{Saglam} to solve the  1D viscous Burgers' equation.
A sixth-order compact finite difference scheme and a Runge-Kutta scheme explicit method were applied by Sari and G$\ddot{\rm o}$rarslan~\cite{Sari} for spatial discretization and time integration respectively for solving the viscous Burgers' equation.
Kutluay et al. \cite{Kutluay1999} presented finite difference solution and analytical solution of the finite difference approximations to the 1D Burgers equation.
The  Hopf-Cole transformation and the standard three-point fourth-order compact finite difference method was applied by Liao~\cite{Liao1,Liao2}  to solve the viscous Burgers' equation.

For the multidimensional coupled Burgers' equations, there are also a lot of studies~\cite{Huang,Bahadir,Liao3,Campos,ChenBuyun.{2017}}.
For instance, a modified local Crank-Nicolson scheme was developed by Huang and Abduwali~\cite{Huang}. A fully implicit finite difference scheme was utilized by Bahadir~\cite{Bahadir} to solve the 2D  viscous Burgers' equations.  Based on the 2D Hopf-Cole transformation, the Crank-Nicolson  scheme together with standard fourth-order explicit compact finite difference method were utilized by Liao~\cite{Liao3} to solve the 2D  viscous Burgers' equations. However, the assumption that the vector consisting of unknown variables  is irrotational, i.e., potential symmetry condition~\cite{Liao3} is needed for the 2D Hopf-Cole transformation, while this condition is only satisfied at certain special cases as mentioned in~\cite{Liao3}. Therefore, the general 2D coupled viscous Burgers' equations can not be solved by the method of Liao~\cite{Liao3}. A Crank-Nicolson scheme together with a linearized approach were developed by Campos et al.~\cite{Campos}. Recently, a linearized high-order accurate CCD-ADI method was proposed in \cite{ChenBuyun.{2017}} to solve 2D/3D coupled Burgers' equations.

The aim of this paper is to derive a  high-order accurate and efficient numerical method for solving the  multi-dimensional coupled Burgers' equations.  To achieve this goal, the third-order total variation diminishing Runge-Kutta (TVD-RK3) scheme and the three-point sixth-order combined compact difference (CCD) scheme will be introduced to solve coupled Burgers' equations. The CCD  scheme was first proposed by Chu and Fan~\cite{PChu} to solve 1D steady convection-diffusion equation.  According to fourier analysis, the CCD scheme possesses better spectral resolution than many other existing high-order schemes~\cite{Mahesh1998}.
 For the compact three-point stencils and high-order accuracy, it has been proved to be very efficient to solve different kinds of differential equations~\cite{HSun,Nihei2003,Sun2,Gao2,DDH,Li,He2017,He20172}.

In this paper, we will develop an  efficient spatial sixth-order linearized method to solve multi-dimensional Burgers' equations. In our numerical scheme, we treat the coupled Burgers' equations in their original form, use the famous TVD-RK3 scheme for the temporal discretization, and adopt the three-point sixth-order CCD scheme for the spatial discretization. Moreover, the unique solvability of the linear CCD system for calculating sixth-order approximations to the first and second derivatives is obtained. When using numerical methods based on Hopt-Cole transformations \cite{Liao1}, some potential symmetry conditions should be imposed, which restricts the application of the method. Our method directly solves the original Burgers' equations, avoiding the use of Hopf-Cole transformation, and can be applied to the problems with  general initial boundary conditions. The proposed CCD-TVD method is third-order accurate in time, sixth-order accurate in space, and is a linearized, explicit scheme, which can be implemented easily.

The rest part of this paper is organized as follows. Section \ref{section1} states the initial boundary problem for the 3D coupled Burgers' equations. Section \ref{section3} presents a CCD-TVD method to solve the 3D coupled nonlinear Burgers'  equations.
We provide a theoretical proof for the unique solvability of the CCD system in section \ref{sec4}.
Section \ref{sec5} gives the numerical results for 1D, 2D and 3D problems, and conclusions are provided in the final section.

\section{Multi-dimensional coupled Burgers' equations }\label{section1}
In this paper, we will present the CCD-TVD method for the following 3D coupled nonlinear viscous Burgers' equations,
\begin{align}
u_{t}+uu_{x}+vu_{y}+wu_{z}&=\frac{1}{Re} (u_{xx}+u_{yy}+u_{zz}),\label{equation1} \\
v_{t}+uv_{x}+vv_{y}+wv_{z}&=\frac{1}{Re} (v_{xx}+v_{yy}+v_{zz}),\label{equation2}\\
w_{t}+uw_{x}+vw_{y}+ww_{z}&=\frac{1}{Re} (w_{xx}+w_{yy}+w_{zz}),\label{equation3}
\end{align}
for $ (x,y,z)\in \Omega=[L_x,R_x]\times[L_y,R_y]\times[L_z,R_z],\ t\in [0, T]$,
with initial conditions
\begin{align}\label{initial}
u(x,y,z,0)=\phi_1 (x,y,z),\quad v(x,y,z,0)=\phi_2 (x,y,z),\quad w(x,y,z,0)=\phi_3 (x,y,z),
\end{align}
and boundary conditions
\begin{align}\label{boundary}
u(x,y,z,t)&=\varphi_1(x,y,z,t),\quad v(x,y,z,t)=\varphi_2(x,y,z,t),\quad w(x,y,z,t)=\varphi_3(x,y,z,t),
\end{align}
for $(x,y,z)\in \partial\Omega,\ 0\leq t\leq T$, where $\phi_1(x,y,z)$,  $\phi_2(x,y,z)$, $\phi_3(x,y,z)$, $\varphi_1(x,y,z,t)$, $\varphi_2(x,y,z,t)$, $\varphi_3(x,y,z,t)$ are given  smooth functions, $L_x,R_x,L_y,R_y,L_z,R_z$ are real constants, $Re$ is the Reynolds number. And we assume that the exact solutions $u,v,w$ are all sufficiently smooth functions. The CCD-TVD method presented in
next section can be applied straightforwardly to 1D and 2D nonlinear Burgers¡¯ equations.

\section{CCD-TVD method}\label{section3}

For a positive integer $N$, let $\Delta t= T/N$, $t_n=n\Delta t$. The time interval $[0,T]$ is covered by $\{t_n| n=0,\cdots,N\}$. Let $f^n$  be the
approximation of $f(x, y,z, n\Delta t)$ for an arbitrary function $f(x, y,z, t)$ at time $t=n\Delta t$. The solution domain $\Omega\times [0,T]$ is covered by a uniform mesh
$\Omega_h=\{(x_i,y_j,z_k,t_n)|x_i=L_x+i\Delta x,y_j=L_y+j\Delta y,z_k=L_z+k\Delta z,t_n=n\tau,i=0,\cdots,M_x, j=0,\cdots,M_y, k=0, \cdots, M_z,  n=0,\cdots, N\}$, where $\Delta x=\frac{R_x-L_x}{M_x}, \Delta y=\frac{R_y-L_y}{M_y},   \Delta z=\frac{R_z-L_z}{M_z}$.

\subsection{Temporal discretization}
It is well known that  the TVD-RK3 method has a CFL coefficient equal to 1, which means that it maintains stability under the same time step restriction as forward Euler
scheme. And it is also pointed that the general four stage fourth order Runge-Kutta methods cannot be TVD without introducing an adjoint operator \cite{TVD.{1998}}.
For time discretization, we adopt the TVD-RK3 scheme for Eqs. (\ref{equation1})-(\ref{equation3}) around $t=t_n$ as follows:
\begin{align}
u^{(1)}&=u^{n}+\Delta t\left(-u^{n}u^{n}_{x}-v^{n}u^{n}_{y}-w^{n}u^{n}_{z}+\frac{1}{Re}( u_{xx}^{n}+ u_{yy}^{n}+u_{zz}^{n})\right),\label{equation4}\\
v^{(1)}&=v^{n}+\Delta t\left(-u^{n}v^{n}_{x}-v^{n}v^{n}_{y}-w^{n}v^{n}_{z}+\frac{1}{Re}( v_{xx}^{n}+ v_{yy}^{n}+v_{zz}^{n})\right),\\
w^{(1)}&=w^{n}+\Delta t\left(-u^{n}w^{n}_{x}-v^{n}w^{n}_{y}-w^{n}w^{n}_{z}+\frac{1}{Re}( w_{xx}^{n}+ w_{yy}^{n}+w_{zz}^{n})\right),\label{equation44}
\end{align}
\begin{align}
u^{(2)}&=\frac{3}{4}u^n + \frac{1}{4}u^{(1)}+\frac{1}{4}\Delta t\left(-u^{(1)}u^{(1)}_{x}-v^{(1)}u^{(1)}_{y}-w^{(1)}u^{(1)}_{z}+\frac{1}{Re}( u_{xx}^{(1)}+ u_{yy}^{(1)}+u_{zz}^{(1)})\right) ,\label{equation5}\\
v^{(2)}&=\frac{3}{4}v^n + \frac{1}{4}v^{(1)}+\frac{1}{4}\Delta t\left(-u^{(1)}v^{(1)}_{x}-v^{(1)}v^{(1)}_{y}-w^{(1)}v^{(1)}_{z}+\frac{1}{Re}( v_{xx}^{(1)}+ v_{yy}^{(1)}+v_{zz}^{(1)})\right) ,\\
w^{(2)}&=\frac{3}{4}w^n + \frac{1}{4}w^{(1)}+\frac{1}{4}\Delta t\left(-u^{(1)}w^{(1)}_{x}-v^{(1)}w^{(1)}_{y}-w^{(1)}w^{(1)}_{z}+\frac{1}{Re}( w_{xx}^{(1)}+ w_{yy}^{(1)}+w_{zz}^{(1)})\right), \label{equation55}
\end{align}
\begin{align}
u^{n+1}&=\frac{1}{3}u^n + \frac{2}{3}u^{(2)}+\frac{2}{3}\Delta t\left(-u^{(2)}u^{(2)}_{x}-v^{(2)}u^{(2)}_{y}-w^{(2)}u^{(2)}_{z}+\frac{1}{Re}( u_{xx}^{(2)}+ u_{yy}^{(2)}+u_{zz}^{(2)})\right) , \label{equation66}\\
v^{n+1}&=\frac{1}{3}v^n + \frac{2}{3}v^{(2)}+\frac{2}{3}\Delta t\left(-u^{(2)}v^{(2)}_{x}-v^{(2)}v^{(2)}_{y}-w^{(2)}v^{(2)}_{z}+\frac{1}{Re}( v_{xx}^{(2)}+ v_{yy}^{(2)}+v_{zz}^{(2)})\right) , \\
w^{n+1}&=\frac{1}{3}w^n + \frac{2}{3}w^{(2)}+\frac{2}{3}\Delta t\left(-u^{(2)}w^{(2)}_{x}-v^{(2)}w^{(2)}_{y}-w^{(2)}w^{(2)}_{z}+\frac{1}{Re}( w_{xx}^{(2)}+ w_{yy}^{(2)}+w_{zz}^{(2)})\right).\label{equation6}
\end{align}
Note that the above discretizations are only for time. And as is verified in \cite{TVD.{1998}} that equations above have the truncation error $O(\Delta t^3)$ .

\subsection{Spatial approximation}

\subsubsection{CCD scheme}
The CCD method was originally proposed to solve 1D steady convection-diffusion equation. Since the TVD-RK3 scheme (\ref{equation4})-(\ref{equation6}) is an explicit scheme, we only need to obtain the
 high-order approximation to the first and second derivatives of $u$, $v$, and $w$. Below we just summarize the CCD method, which will be used to compute the
 sixth-order approximations to these derivatives.

Let $u$ be defined on the intervals $[a,b]$. Then we divide the internal $[a,b]$ with mesh size $h=\frac{b-a}{M}$ and $x_{i}=a+ih (i=0,1,\cdots,M)$.
Let $u_i^{'}$ and $u_i^{''}$ be the numerical approximations of $u_x(x_i)$ and $u_{xx}(x_i)$, respectively.
The first and second derivatives $u_i', u_i''$ for all mesh points can be computed as follows.

For internal points $x_i(i=1,2,\cdots,M-1)$, we have
 \begin{align}
\frac{7}{16}(u'_{i+1}+u'_{i-1})+u'_i-\frac{h}{16}(u''_{i+1}-u''_{i-1})-\frac{15}{16h}(u_{i+1}-u_{i-1})=0, \label{equation18} \\
\frac{9}{8h}(u'_{i+1}-u'_{i-1})+u''_i-\frac{1}{8}(u''_{i+1}+u''_{i-1})-\frac{3}{h^2}(u_{i+1}-2u_i+u_{i-1})=0.\label{equation19}
\end{align}
The above two relations are valid up to $O(h^6)$ when the exact solution $u(x)$ is 8-times continuously differentiable.

For boundary points, the three-point one-side fifth-order schemes are used,
 \begin{align}
14u'_0 + 16u'_1 + 2hu''_0 - 4h u''_1 +\frac{1}{h}(31u_0 -32u_1+ u_2) &=0,\label{equation20} \\
14u'_{M} + 16u'_{M-1} -2hu''_{M} + 4h u''_{M-1} -\frac{1}{h}(31u_{M} -32u_{M-1}+ u_{M-2}) &=0.\label{equation21}
\end{align}

Since there are only $2M$ equations above, we also need two additional boundary CCD equations
    \begin{align}
    u^{'}_{0}+2u^{'}_{1}-hu^{''}_{1}\mathbf{+}\frac{1}{2h}(7u_{0}-8u_{1}+u_{2})&=0, \\
    u^{'}_{M}+2u^{'}_{M-1}\mathbf{+}hu^{''}_{M-1}-\frac{1}{2h}(7u_{M}-8u_{M-1}+u_{M-2})&=0, \label{abc}
    \end{align}
    to compute the all $2(M+1)$ unknown derivatives $u^{'}_{i}$,$u^{''}_{i}(i=0,1,\cdots,M)$.

 In order to write the CCD system (\ref{equation18})-(\ref{abc}) in a concise style, we order all the unknowns in the natural column-wise sense,
\begin{equation}\label{vv}
\mathbf{v} = \left[ \begin{array}{l}
\mathbf{u}_x \\
\mathbf{u}_{xx}
\end{array}
\right],
\end{equation}
where the unknown derivative vectors {are}
\begin{equation*}
  \mathbf{u_x} = [u_0^{'},u_1^{'},\cdots,u_M^{'}]^T, \; \mathbf{u_{xx}} = [u_0^{''},u_1^{''},\cdots,u_M^{''}]^T,
\end{equation*}
and the known function vector is
\begin{equation*}
  \mathbf{u} = [u_0,u_1,\cdots,u_M]^T.
\end{equation*}

Let $m\triangleq M+1$, then we can rewrite the CCD system (\ref{equation18})-(\ref{abc}) as the following $2\times 2$ block linear system,
\begin{equation}\label{sys}
A\mathbf{v} = \left[
    \begin{array}{cc}
    A1&A2\\
    A3&A4\\
    \end{array}\right]
\cdot
\left[ \begin{array}{c}
\mathbf{u}_x \\
\mathbf{u}_{xx}
\end{array}
\right]
=
\left[
    \begin{array}{c}
    B1\\
    B2\\
    \end{array}\right]\cdot \mathbf{u},
\end{equation}
    where
\begin{align}
  A1&=\left[ \begin{matrix}14&16&0&\cdots&0\\ \frac{7}{16}&1&\frac{7}{16}&\quad&\vdots\\ 0&\ddots&\ddots&\ddots&0\\ \vdots&\quad&\frac{7}{16}&1&\frac{7}{16}\\ 0&\cdots&0&16&14\\  \end{matrix}\right]_{m\times m}, \quad A2=\left[ \begin{matrix}2h&-4h&0&\cdots&0\\ \frac{h}{16}&0&-\frac{h}{16}&\quad&\vdots\\ 0&\ddots&\ddots&\ddots&0\\ \vdots&\quad&\frac{h}{16}&0&-\frac{h}{16}\\ 0&\cdots&0&{4h}&{-2h}\\  \end{matrix}\right]_{m\times m},\\
A3&=\left[ \begin{matrix}1&2&0&\cdots&0\\ {-\frac{9}{8h}}&0&\frac{9}{8h}&\quad&\vdots\\ 0&\ddots&\ddots&\ddots&0\\ \vdots&\quad&{-\frac{9}{8h}}&0&\frac{9}{8h}\\ 0&\cdots&0&2&1\\  \end{matrix}\right]_{m\times m}, \quad
A4=\left[ \begin{matrix}0&-h&0&\cdots&0\\ {-\frac{1}{8}}&1&-\frac{1}{8}&\quad&\vdots\\ 0&\ddots&\ddots&\ddots&0\\ \vdots&\quad&{-\frac{1}{8}}&1&-\frac{1}{8}\\ 0&\cdots&0&{h}&0\\  \end{matrix}\right]_{m\times m},
\end{align}
and
\begin{align}
  B1=\left[ \begin{matrix}{-\frac{31}{h}}&{\frac{32}{h}}&{-\frac{1}{h}}&\cdots&0\\ -\frac{5}{16h}&0&\frac{5}{16h}&\quad&\vdots\\ 0&\ddots&\ddots&\ddots&0\\ \vdots&\quad&-\frac{5}{16h}&0&\frac{5}{16h}\\ 0&\cdots&\frac{1}{h}&-\frac{32}{h}&\frac{31}{h}\\  \end{matrix}\right]_{m\times m},\quad B2=\left[ \begin{matrix}{-\frac{7}{2h}}&{\frac{8}{2h}}&{-\frac{1}{2h}}&\cdots&0\\ {\frac{3}{h^{2}}}&{-\frac{6}{h^{2}}}&\frac{3}{h^{2}}&\quad&\vdots\\ 0&\ddots&\ddots&\ddots&0\\ \vdots&\quad&{\frac{3}{h^{2}}}&{-\frac{6}{h^{2}}}&\frac{3}{h^{2}}\\ 0&\cdots&\frac{1}{2h}&-\frac{8}{2h}&\frac{7}{2h}\\  \end{matrix}\right]_{m\times m}.
\end{align}

    Denote $IAB=\left[
    \begin{array}{cc}
    A1&A2\\
    A3&A4\\
    \end{array}\right]^{-1}\left[
    \begin{array}{c}
    B1\\
    B2\\
    \end{array}\right]$. Then, we can get $\mathbf{u_x},  \mathbf{u_{xx}}$ immediately by performing one matrix-vector multiplication $IAB\cdot\mathbf{u}$. Here $A1, A2, A3, A4$ and $B1, B2$ are all constant tridiagonal matrix independent of time level, thus the matrix $IAB$ only need be calculated once.

    The above CCD scheme (\ref{equation18})-(\ref{abc}) is sixth-order accurate~\cite{PChu,Li}. Combining the third-order accuracy in time, the overall accuracy of $u$, $v$ and $w$ is $O ((\Delta t)^3 + (\Delta x)^6 + (\Delta y)^6+ (\Delta z)^6)$.


\subsubsection{Implementation of  CCD-TVD scheme}

Algorithm \ref{alg1} describes the detailed implementation of the CCD-ADI method (\ref{equation4})-(\ref{equation6}) together with (\ref{equation18})-(\ref{abc}) for solving the 3D coupled nonlinear Burgers¡¯ equations (\ref{equation1})-(\ref{boundary}).

\begin{algorithm}[!tbp]         
\caption{ CCD-TVD algorithm for the 3D coupled nonlinear Burgers' equations}             
\label{alg1}
\begin{algorithmic}[1]
\STATE Compute $u^0, v^0$ and $w^0$ for all mesh points by using~(\ref{initial}).
 \FOR {$n=0$ to $N-1$}
    \STATE Step 1:  Compute the right-hand side of (\ref{equation4})-(\ref{equation44}) for all mesh
points.\\
 \begin{description}
  \item[\quad a)] Do $i=0,\cdots,M_x$, $j=0,\cdots,M_y$
compute $(u^{n}_z)_{i,j,k}$,$(u^{n}_{zz})_{i,j,k}$,$(v^{n}_z)_{i,j,k}$, $(v^{n}_{zz})_{i,j,k}$, $(w^{n}_z)_{i,j,k}$, $(w^{n}_{zz})_{i,j,k}$ ($k=0,\cdots,M_z$) by the CCD scheme (\ref{equation18})-(\ref{abc});
  \item[\quad b)] Do $i=0,\cdots,M_x$, $k=0,\cdots,M_z$
compute $(u^{n}_y)_{i,j,k}$,$(u^{n}_{yy})_{i,j,k}$, $(v^{n}_y)_{i,j,k}$, $(v^{n}_{yy})_{i,j,k}$, $(w^{n}_y)_{i,j,k}$, $(w^{n}_{yy})_{i,j,k}$ ($j=0,\cdots,M_y$) by the CCD scheme (\ref{equation18})-(\ref{abc});
  \item[\quad c)] Do $j=0,\cdots,M_y$, $k=0,\cdots,M_z$
compute $(u^{n}_x)_{i,j,k}$,$(u^{n}_{xx})_{i,j,k}$,$(v^{n}_x)_{i,j,k}$,$(v^{n}_{xx})_{i,j,k}$,$(w^{n}_x)_{i,j,k}$,$(w^{n}_{xx})_{i,j,k}$($i=0,\cdots,M_x$) by the CCD scheme (\ref{equation18})-(\ref{abc});
   \item[\quad d)] Compute $u^{(1)}$, $v^{(1)}$ and $w^{(1)}$ using (\ref{equation4})-(\ref{equation44}) for all mesh points.
\end{description}

 \STATE  Step 2: Compute the right-hand side of (\ref{equation5})-(\ref{equation55}) for all mesh
points.
\begin{description}
  \item[\quad a)] Do $i=0,\cdots,M_x$, $j=0,\cdots,M_y$
compute $(u^{(1)}_z)_{i,j,k}$, $(u^{(1)}_{zz})_{i,j,k}$, $(v^{(1)}_z)_{i,j,k}$, $(v^{(1)}_{zz})_{i,j,k}$, $(w^{(1)}_z)_{i,j,k}$, $(w^{(1)}_{zz})_{i,j,k}$($k=0,\cdots,M_z$) by the CCD scheme (\ref{equation18})-(\ref{abc});
  \item[\quad b)] Do $i=0,\cdots,M_x$, $k=0,\cdots,M_z$
compute $(u^{(1)}_y)_{i,j,k}$, $(u^{(1)}_{yy})_{i,j,k}$, $(v^{(1)}_y)_{i,j,k}$, $(v^{(1)}_{yy})_{i,j,k}$, $(w^{(1)}_y)_{i,j,k}$, $(w^{(1)}_{yy})_{i,j,k}$($j=0,\cdots,M_y$) by the CCD scheme (\ref{equation18})-(\ref{abc});
  \item[\quad c)] Do $j=0,\cdots,M_y$, $k=0,\cdots,M_z$
compute $(u^{(1)}_x)_{i,j,k}$, $(u^{(1)}_{xx})_{i,j,k}$, $(v^{(1)}_x)_{i,j,k}$, $(v^{(1)}_{xx})_{i,j,k}$, $(w^{(1)}_x)_{i,j,k}$, $(w^{(1)}_{xx})_{i,j,k}$($i=0,\cdots,M_x$) by the CCD scheme (\ref{equation18})-(\ref{abc});
   \item[\quad d)] Compute $u^{(2)}$, $v^{(2)}$ and $w^{(2)}$ using  (\ref{equation5})-(\ref{equation55})  for all mesh points.
\end{description}

\STATE Step 3: Compute the right-hand side of (\ref{equation66})-(\ref{equation6}) for all mesh
points.\\
\begin{description}
  \item[\quad a)] Do $i=0,\cdots,M_x$, $j=0,\cdots,M_y$
compute $(u^{(2)}_z)_{i,j,k}$, $(u^{(2)}_{zz})_{i,j,k}$, $(v^{(2)}_z)_{i,j,k}$, $(v^{(2)}_{zz})_{i,j,k}$, $(w^{(2)}_z)_{i,j,k}$, $(w^{(2)}_{zz})_{i,j,k}$($k=0,\cdots,M_z$) by the CCD scheme (\ref{equation18})-(\ref{abc});
  \item[\quad b)] Do $i=0,\cdots,M_x$, $k=0,\cdots,M_z$
compute $(u^{(2)}_y)_{i,j,k}$, $(u^{(2)}_{yy})_{i,j,k}$, $(v^{(2)}_y)_{i,j,k}$, $(v^{(2)}_{yy})_{i,j,k}$, $(w^{(2)}_y)_{i,j,k}$, $(w^{(2)}_{yy})_{i,j,k}$($j=0,\cdots,M_y$) by the CCD scheme (\ref{equation18})-(\ref{abc});
  \item[\quad c)] Do $j=0,\cdots,M_y$, $k=0,\cdots,M_z$
compute $(u^{(2)}_x)_{i,j,k}$, $(u^{(2)}_{xx})_{i,j,k}$, $(v^{(2)}_x)_{i,j,k}$, $(v^{(2)}_{xx})_{i,j,k}$, $(w^{(2)}_x)_{i,j,k}$, $(w^{(2)}_{xx})_{i,j,k}$($i=0,\cdots,M_x$) by the CCD scheme (\ref{equation18})-(\ref{abc});
   \item[\quad d)] Compute $u^{n+1}$, $v^{n+1}$ and $w^{n+1}$ using (\ref{equation66})-(\ref{equation6}) for all mesh points.

\end{description}
\STATE Step 4: Update the boundary values of $u^{n+1}$, $v^{n+1}$ and $w^{n+1}$ by using~(\ref{boundary}).

 \ENDFOR
\end{algorithmic}

\end{algorithm}

\begin{remark}
  The CCD-TVD method presented in the following section can be applied  straightforwardly to 2D coupled nonlinear Burgers' equations and 1D Burgers' equation.
\end{remark}

\begin{remark}
 The CCD-TVD method doesn't  rely on the condition that the initial vorticity is zero.  Thus, it can be applied to the multi-dimensional coupled Burgers' equations with general initial and boundary conditions.
\end{remark}

\section{Unique solvability of the CCD system}\label{sec4}
    In this section,we provide a theoretical proof for the uniqueness of the numerical solution of (\ref{sys}).

\begin{definition}\label{def1}
  $C^{3}_{m}(a,b,c,d,e,f,g)$ is an m-by-m semi-circulant matrix determined by 3 circulant elements a,b,c and 4 boundary elements d,e,f and g, which is defined as follows,\\
\[\centering
 C^{3}_{m}(a,b,c,d,e,f,g)=\left[ \begin{matrix}d&e&0&\cdots&0\\ a&b&c&\quad&\vdots\\ 0&\ddots&\ddots&\ddots&0\\ \vdots&\quad&a&b&c\\ 0&\cdots&0&g&f\\  \end{matrix}\right]_{m \times m} \in R^{mxm}.\]
 \end{definition}

\begin{lemma}\label{lem1}
  \begin{equation}
 C^{3}_{m}(a_1,a_2,a_3,a_4,a_5,a_6,a_7)C^{3}_{m}(b_1,b_2,b_3,b_4,b_5,b_6,b_7)=C^{5}_{m}(c_1,c_2,\cdots,c_5,d_1,d_2,\cdots,d_8,e_1,e_2,\cdots,e_6),
 \end{equation}
\end{lemma}
 where
\[
 C^{5}_{m}(c_1,c_2,\cdots,c_5,d_1,d_2,\cdots,d_8,e_1,e_2,\cdots,e_6)=\left[ \begin{matrix}e_1&e_2&e_3&0&\cdots&\cdots&\cdots&0\\ d_1&d_2&d_3&d_4&0&\cdots&\cdots&0\\ c_1&c_2&c_3&c_4&c_5&0&\cdots&0\\ 0&\ddots&\ddots&\ddots&\ddots&\ddots&\ddots&\vdots\\ \vdots&\ddots&\ddots&\ddots&\ddots&\ddots&\ddots&0\\ 0&\cdots&0&c_1&c_2&c_3&c_4&c_5 \\ 0&\cdots&\cdots&0&d_5&d_6&d_7&d_8 \\ 0&\cdots&\cdots&\cdots&0&e_4&e_5&e_6 \end{matrix}\right]_{m \times m} \in R^{mxm},\]
 and
 \[\begin{split}
    \centering
     e_1 &=a_4b_4+a_5b_1,\quad
     e_2 =a_4b_5+a_5b_2,\quad
     e_3 =a_5b_3,\\
     e_4 &=a_7b_1,\quad
     e_5 =a_6b_7+a_7b_2,\quad
     e_6 =a_7b_3+a_6b_6,\\
     d_1 &=a_1b_4+a_2b_1,\quad
     d_2  =a_1b_5+a_2b_2+a_3b_1,\quad
     d_3  =a_2b_3+a_3b_2,\quad
     d_4  =a_3b_3,\\
     d_5 &=a_1b_1,\quad
     d_6 =a_1b_2+a_2b_1,\quad
     d_7 =a_1b_3+a_2b_2+a_3b_7,\quad
     d_8 =a_2b_3+a_3b_6,\\
     c_1 &=a_1b_1,\quad
     c_2 =a_1b_2+a_2b_1,\quad
     c_3 =a_1b_3+a_2b_2+a_3b_1,\\
     c_4 &=a_2b_3+a_3b_2,\quad
     c_5 =a_3b_3.
\end{split}
 \]
\noindent {\bf Proof}. This lemma can be verified through direct computation.\qed

\begin{lemma}\label{lem2}
Assume that $A,B,C,D$ are $n\times n$ matrices.
If $AC=CA$, then
$$\left|\begin{array}{cc}
      A & B  \\
       C & D
     \end{array}
\right|=\big|AD-CB\big|.$$
\end{lemma}
\noindent {\bf Proof}. This lemma can be found in any elemental Linear algebra text book. \qed

\begin{theorem}
  The CCD linear system (\ref{sys}) has a unique solution.
\end{theorem}
\noindent {\bf Proof}. Using Definition \ref{def1}, from (\ref{sys}) we know
\begin{equation}
\left[
    \begin{array}{cc}
    A1&A2\\
    A3&A4\\
    \end{array}\right]
\cdot
\left[ \begin{array}{c}
\mathbf{u}_x \\
\mathbf{u}_{xx}
\end{array}
\right]
=
\left[
    \begin{array}{c}
    B1\\
    B2\\
    \end{array}\right]\cdot \mathbf{u},
\end{equation}
 where \\
 \begin{align*}
 A1&=C^{3}_{m}(\frac{7}{16},1,\frac{7}{16},14,16,14,16),\;\;A2=C^{3}_{m}(\frac{h}{16},0,-\frac{h}{16},2h,-4h,-2h,4h),\\
 A3&=C^{3}_{m}(-\frac{9}{8h},0,\frac{9}{8h},1,2,1,2),\;\;A4=C^{3}_{m}(-\frac{1}{8},1,-\frac{1}{8},0,-h,0,h).
 \end{align*}
 In the following, we will show that the coefficient matrix of the above $2\times 2$ block linear system is nonsingular.

The determinant of the coefficient matrix can be calculated as follows:
 {\small
 \begin{align*}
 &\left|\begin{array}{ll}
 C^{3}_{m}(\frac{7}{16},1,\frac{7}{16},14,16,14,16) & C^{3}_{m}(\frac{h}{16},0,-\frac{h}{16},2h,-4h,-2h,4h)\\
  C^{3}_{m}(-\frac{9}{8h},0,\frac{9}{8h},1,2,1,2)   & C^{3}_{m}(-\frac{1}{8},1,-\frac{1}{8},0,-h,0,h)\\
\end{array}\right|\\
=&(\frac{1}{h})^{m}\left|\begin{array}{ll}
 C^{3}_{m}(\frac{7}{16},1,\frac{7}{16},14,16,14,16) & C^{3}_{m}(\frac{h}{16},0,-\frac{h}{16},2h,-4h,-2h,4h)\\
          C^{3}_{m}(-\frac{9}{8},0,\frac{9}{8},h,2h,h,2h)   & C^{3}_{m}(-\frac{h}{8},h,-\frac{h}{8},0,-h^2,0,h^2)\\
\end{array}\right|\\
=&(\frac{1}{h})^{m}h^m\left|\begin{array}{ll}
 C^{3}_{m}(\frac{7}{16},1,\frac{7}{16},14,16,14,16) & C^{3}_{m}(\frac{1}{16},0,-\frac{1}{16},2,-4,-2,4)\\
          C^{3}_{m}(-\frac{9}{8},0,\frac{9}{8},h,2h,h,2h)   & C^{3}_{m}(-\frac{1}{8},1,-\frac{1}{8},0,-h,0,h)\\
\end{array}\right|\\
=&h^2\left|\begin{array}{ll}
 C^{3}_{m}(\frac{7}{16},1,\frac{7}{16},14,16,14,16) & C^{3}_{m}(\frac{1}{16},0,-\frac{1}{16},2,-4,-2,4)\\
          C^{3}_{m}(-\frac{9}{8},0,\frac{9}{8},1,2,1,2)   & C^{3}_{m}(-\frac{1}{8},1,-\frac{1}{8},0,-1,0,1)\\
\end{array}\right|\\
=&h^2\left|\begin{array}{ll}
 C^{3}_{m}(\frac{7}{16},1,\frac{7}{16},14,16,14,16) & C^{3}_{m}(\frac{1}{16},0,-\frac{1}{16},2,-4,-2,4)\\
          C^{3}_{m}(-\frac{9}{8},0,\frac{9}{8},1,2,1,2)+\frac{6}{7}\sqrt{7}C^{3}_{m}(\frac{7}{16},1,\frac{7}{16},14,16,14,16)& C^{3}_{m}(-\frac{1}{8},1,-\frac{1}{8},0,-1,0,1)+\frac{6}{7}\sqrt{7}C^{3}_{m}(\frac{1}{16},0,-\frac{1}{16},2,-4,-2,4)\\
\end{array}\right|\\
=&h^2\left|\begin{array}{ll}
 C^{3}_{m}(\frac{7}{16},1,\frac{7}{16},14,16,14,16) & C^{3}_{m}(\frac{1}{16},0,-\frac{1}{16},2,-4,-2,4)\\
          C^{3}_{m}(\frac{-9+3\sqrt{7}}{8},\frac{6\sqrt{7}}{7},\frac{3\sqrt{7}+9}{8},1+12\sqrt{7},\frac{14+96\sqrt{7}}{7},1+12\sqrt{7},\frac{14+96\sqrt{7}}{7})& C^{3}_{m}(\frac{3\sqrt{7}-7}{56},1,-\frac{3\sqrt{7}+7}{56},\frac{12\sqrt{7}}{7},-\frac{24\sqrt{7}+7}{7},-\frac{12\sqrt{7}}{7},\frac{24\sqrt{7}+7}{7})\\
\end{array}\right|\\
=&h^2\left|\begin{array}{ll}
 C^{3}_{m}(\frac{7}{16},1,\frac{7}{16},14,16,14,16)+3\sqrt{7}C^{3}_{m}(\frac{1}{16},0,-\frac{1}{16},2,-4,-2,4) & C^{3}_{m}(\frac{1}{16},0,-\frac{1}{16},2,-4,-2,4)\\
          C^{3}_{m}(\frac{-9+3\sqrt{7}}{8},\frac{6\sqrt{7}}{7},\frac{3\sqrt{7}+9}{8},1+12\sqrt{7},\frac{14+96\sqrt{7}}{7},1+12\sqrt{7},\frac{14+96\sqrt{7}}{7})& C^{3}_{m}(\frac{3\sqrt{7}-7}{56},1,-\frac{3\sqrt{7}+7}{56},\frac{12\sqrt{7}}{7},-\frac{24\sqrt{7}+7}{7},-\frac{12\sqrt{7}}{7},\frac{24\sqrt{7}+7}{7})\\
          +3\sqrt{7}C^{3}_{m}(\frac{3\sqrt{7}-7}{56},1,-\frac{3\sqrt{7}+7}{56},\frac{12\sqrt{7}}{7},-\frac{24\sqrt{7}+7}{7},-\frac{12\sqrt{7}}{7},\frac{24\sqrt{7}+7}{7})& \\
\end{array}\right|\\
=&h^2\left|\begin{array}{ll}
 C^{3}_{m}(\frac{7+3\sqrt{7}}{16},1,\frac{7-3\sqrt{7}}{16},14+6\sqrt{7},16-12\sqrt{7},14-6\sqrt{7},16+12\sqrt{7}) & C^{3}_{m}(\frac{1}{16},0,-\frac{1}{16},2,-4,-2,4)\\
          C^{3}_{m}(0,\frac{27\sqrt{7}}{7},0,37+12\sqrt{7},\frac{75\sqrt{7}-490}{7},-35+12\sqrt{7},\frac{117\sqrt{7}+518}{7})& C^{3}_{m}(\frac{3\sqrt{7}-7}{56},1,-\frac{3\sqrt{7}+7}{56},\frac{12\sqrt{7}}{7},-\frac{24\sqrt{7}+7}{7},-\frac{12\sqrt{7}}{7},\frac{24\sqrt{7}+7}{7})\\
\end{array}\right|\\
=&h^2(\frac{27}{\sqrt{7}})^{m-2}\left|\begin{array}{ll}
 C^{3}_{m}(\frac{7+3\sqrt{7}}{16},1,\frac{7-3\sqrt{7}}{16},14+6\sqrt{7},16-12\sqrt{7},14-6\sqrt{7},16+12\sqrt{7}) & C^{3}_{m}(\frac{1}{16},0,-\frac{1}{16},2,-4,-2,4)\\
          C^{3}_{m}(0,1,0,37+12\sqrt{7},\frac{75\sqrt{7}-490}{7},-35+12\sqrt{7},\frac{117\sqrt{7}+518}{7})      &    C^{3}_{m}(\frac{21-7\sqrt{7}}{1512},\frac{\sqrt{7}}{27},-\frac{21+7\sqrt{7}}{1512},\frac{12\sqrt{7}}{7},-\frac{24\sqrt{7}+7}{7},-\frac{12\sqrt{7}}{7},\frac{24\sqrt{7}+7}{7})\\
\end{array}\right|\\
=&h^2(\frac{27}{\sqrt{7}})^{m-2}\left|\begin{array}{ll}
 C^{3}_{m}(\frac{7+3\sqrt{7}}{16},1,\frac{7-3\sqrt{7}}{16},14+6\sqrt{7},16-12\sqrt{7},14-6\sqrt{7},16+12\sqrt{7}) & C^{3}_{m}(\frac{1}{16},0,-\frac{1}{16},2,-4,-2,4)\\
C^{3}_{m}(0,1,0,37+12\sqrt{7},\frac{75\sqrt{7}-490}{7},-35+12\sqrt{7},\frac{117\sqrt{7}+518}{7}) & C^{3}_{m}(\frac{21-7\sqrt{7}}{1512},\frac{\sqrt{7}}{27},-\frac{21+7\sqrt{7}}{1512},\frac{12\sqrt{7}}{7},-\frac{24\sqrt{7}+7}{7},-\frac{12\sqrt{7}}{7},\frac{24\sqrt{7}+7}{7})\\
-\frac{-5(-98+15\sqrt{7})}{28(-4+3\sqrt{7}}C^{3}_{m}(0,0,0,14+6\sqrt{7},16-12\sqrt{7},0,0) & -\frac{-5(-98+15\sqrt{7})}{28(-4+3\sqrt{7}}C^{3}_{m}(0,0,0,2,-4,0,0)\\
-\frac{(518+117\sqrt{7})}{28(4+3\sqrt{7})}C^{3}_{m}(0,0,0,0,0,14-6\sqrt{7},16+12\sqrt{7})     &    -\frac{(518+117\sqrt{7})}{28(4+3\sqrt{7})}C^{3}_{m}(0,0,0,0,0,-2,4)\\
\end{array}\right|\\
=&h^2(\frac{27}{\sqrt{7}})^{m-2}\left|\begin{array}{ll}
 C^{3}_{m}(\frac{7+3\sqrt{7}}{16},1,\frac{7-3\sqrt{7}}{16},14+6\sqrt{7},16-12\sqrt{7},14-6\sqrt{7},16+12\sqrt{7}) & C^{3}_{m}(\frac{1}{16},0,-\frac{1}{16},2,-4,-2,4)\\
          C^{3}_{m}(0,1,0,\frac{-57-9\sqrt{7}}{-8+6\sqrt{7}},0,\frac{57-9\sqrt{7}}{8+6\sqrt{7}},0)      &    C^{3}_{m}(\frac{21-7\sqrt{7}}{1512},\frac{\sqrt{7}}{27},-\frac{21+7\sqrt{7}}{1512},\frac{2-3\sqrt{7}}{-8+6\sqrt{7}},\frac{2}{-4+3\sqrt{7}},\frac{2+3\sqrt{7}}{8+6\sqrt{7}},\frac{2}{4+3\sqrt{7}})\\
\end{array}\right|\\
=&C_0 h^2\left|\begin{array}{ll}
 C^{3}_{m}(\frac{7+3\sqrt{7}}{16},1,\frac{7-3\sqrt{7}}{16},14+6\sqrt{7},16-12\sqrt{7},14-6\sqrt{7},16+12\sqrt{7}) & C^{3}_{m}(\frac{1}{16},0,-\frac{1}{16},2,-4,-2,4)\\
          C^{3}_{m}(0,1,0,1,0,1,0)      &    C^{3}_{m}(\frac{21-7\sqrt{7}}{1512},\frac{\sqrt{7}}{27},-\frac{21+7\sqrt{7}}{1512},\frac{-2+3\sqrt{7}}{57+9\sqrt{7}},-\frac{4}{57+9\sqrt{7}},-\frac{2+3\sqrt{7}}{-57+9\sqrt{7}},-\frac{4}{-57+9\sqrt{7}})\\
\end{array}\right|\\
\end{align*}
}
where
\begin{equation}
  C_0 = (\frac{27}{\sqrt{7}})^{m-2}(\frac{-57-9\sqrt{7}}{-8+6\sqrt{7}})(\frac{57-9\sqrt{7}}{8+6\sqrt{7}}),
\end{equation}
and
$C^{3}_{m}(0,1,0,1,0,1,0)$ is identity matrix which is commutative to any matrix of the same size.

Using Lemma \ref{lem1} and Lemma \ref{lem2}, we have
{\small \begin{align*}
& \left|\begin{array}{ll}
 C^{3}_{m}(\frac{7+3\sqrt{7}}{16},1,\frac{7-3\sqrt{7}}{16},14+6\sqrt{7},16-12\sqrt{7},14-6\sqrt{7},16+12\sqrt{7}) & C^{3}_{m}(\frac{1}{16},0,-\frac{1}{16},2,-4,-2,4)\\
          C^{3}_{m}(0,1,0,1,0,1,0)      &    C^{3}_{m}(\frac{21-7\sqrt{7}}{1512},\frac{\sqrt{7}}{27},-\frac{21+7\sqrt{7}}{1512},\frac{-2+3\sqrt{7}}{57+9\sqrt{7}},-\frac{4}{57+9\sqrt{7}},-\frac{2+3\sqrt{7}}{-57+9\sqrt{7}},-\frac{4}{-57+9\sqrt{7}})\\
\end{array}\right|\\
=& \left| \begin{array}{rrrrrrr}
   \frac{34\sqrt{7}+33}{513+81\sqrt{7}}& \frac{288+160\sqrt{7}}{513+81\sqrt{7}}&\frac{(-4+3\sqrt{7})(7\sqrt{7}+21)}{378} &        0          &   0  &   \cdots   &    0    \\
 \frac{17\sqrt{7}}{4104+648\sqrt{7}}&\frac{567+695\sqrt{7}}{32832+5184\sqrt{7}}&\frac{5\sqrt{7}}{432}& \frac{\sqrt{7}}{1728}  &   0  &  \cdots  &    0     \\
\frac{\sqrt{7}}{1728}& \frac{5\sqrt{7}}{432}&\frac{\sqrt{7}}{36}&  \frac{5\sqrt{7}}{432}&    \frac{\sqrt{7}}{1728} & \cdots & 0     \\
\vdots& \ddots &  \ddots     &   \ddots   & \ddots &  \ddots &    \vdots     \\
     0  &\cdots & \frac{\sqrt{7}}{1728}&  \frac{5\sqrt{7}}{432}& \frac{\sqrt{7}}{36}&  \frac{5\sqrt{7}}{432}&     \frac{\sqrt{7}}{1728}  \\
 0  &\cdots& 0&\frac{\sqrt{7}}{1728}& \frac{5\sqrt{7}}{432}& \frac{-695\sqrt{7}+567}{-32832+5184\sqrt{7}}& \frac{-17\sqrt{7}}{-4104+648\sqrt{7}} \\
     0      &  \cdots      &     0      &   0         &    \frac{(4+3\sqrt{7})(-7\sqrt{7}+21)}{378}&\frac{-160\sqrt{7}+288}{-513+81\sqrt{7}}&\frac{-34\sqrt{7}+33}{-513+81\sqrt{7}}  \\
\end{array} \right|.
 \end{align*}}

 Next, we will apply a series of elementary transformation to the above five-diagonal matrix to show that it is nonsingular.
The detailed procedures are described as follows:
\begin{enumerate}[Step 1:]
  \item $$r_1 \leftarrow r_1-\frac{1459440\sqrt{7}+8541848}{598633} r_2 +\frac{94986\sqrt{7}+1563660}{598633}
 r_3-\frac{55035\sqrt{7}+ 2841419}{3591798} r_4, $$
  \item $$c_3 \leftarrow c_3 + m_{31} c_1,$$
  where $$m_{31}=\frac{161433961059948782743125\sqrt{7}}{638365996543160612814848} + \frac{386982051292812235294125}{319182998271580306407424},$$
  \item $$c_5 \leftarrow c_5 + m_{51} c_1,$$
  where $$m_{51} = \frac{128698051973330045562453\sqrt{7}}{638365996543160612814848} + \frac{320818233644731054212525}{319182998271580306407424},$$
  \item $$r_2 \leftarrow r_2- m_{23} r_3, $$
  where $$m_{23}=\frac{455021090726735024960954900487104822899500\sqrt{7}}{57625186587725827104855703164883727826611651} + \frac{98624527354971701012117024209404389139084220}{172875559763177481314567109494651183479834953},$$
  \item $$c_{n-1} \leftarrow c_{n-1} -  \frac{3}{2}c_n,$$
  \item $$c_{n-2} \leftarrow c_{n-2} -  \frac{2269-570\sqrt{7}}{1490}c_n,$$
  \item $$ r_{n-1} \leftarrow r_{n-1}- \frac{1}{10}r_{n-2} - \frac{1}{10} r_n.$$
\end{enumerate}

Finally, we obtain
\begin{align}
&\left| \begin{array}{rrrrrrrrr}
    T_1    &                T_2     &           0 &  T_3    &         0&       T_4&         0&    \cdots&         0\\
    T_5    &                T_6     &           0 &  T_7    &       T_8&         0&         0&    \cdots&         0\\
\frac{\sqrt{7}}{1728}&\frac{5\sqrt{7}}{432}&        T_9 &\frac{5\sqrt{7}}{432}&      T_{10}&         0&         0&    \cdots&         0\\
     0 &  \frac{\sqrt{7}}{1728} &\frac{5\sqrt{7}}{432}&\frac{\sqrt{7}}{36}&\frac{5\sqrt{7}}{432}&  \frac{\sqrt{7}}{1728} &         0&    \cdots&         0\\
    \vdots &                  \quad &      \ddots &   \ddots&    \ddots&    \ddots&    \ddots&    \quad &    \vdots\\
         0 &                  \cdots&           0 &  \frac{\sqrt{7}}{1728} &\frac{5\sqrt{7}}{432}&\frac{\sqrt{7}}{36}&\frac{5\sqrt{7}}{432}&  \frac{\sqrt{7}}{1728} &         0\\
         0 &                  \cdots&           0 &        0&\frac{\sqrt{7}}{1728}&\frac{5\sqrt{7}}{432}&      T_{11}&      T_{12}&\frac{\sqrt{7}}{1728}\\
         0 &                  \cdots&           0 &        0&      T_{13}&  -\frac{\sqrt{7}}{1728} &    T_{14}&      T_{15}&      T_{16}\\
         0 &                  \cdots&           0 &        0&         0&         0&      T_{17}&      T_{18}&      T_{19}\\
\end{array} \right|\label{aaa}\\
\approx &\left| \begin{array}{rrrrrrrrr}
    0.0135 &  -0.0013 &        0 &        0&         0&   -0.0013&         0&    \cdots&         0\\
    0.0068 &   0.0336 &        0 &  -0.0166&    0.0096&         0&         0&    \cdots&         0\\
    0.0015 &   0.0306 &   0.0764 &   0.0306&    0.0039&         0&         0&    \cdots&         0\\
         0 &   0.0015 &   0.0306 &   0.0735&    0.0306&    0.0015&         0&    \cdots&         0\\
    \vdots &    \quad &   \ddots &   \ddots&    \ddots&    \ddots&    \ddots&    \quad &    \vdots\\
         0 &    \cdots&        0 &   0.0015&    0.0306&    0.0735&    0.0306&    0.0015&         0\\
         0 &    \cdots&        0 &        0&    0.0015&    0.0306&    0.0727&    0.0283&    0.0015\\
         0 &    \cdots&        0 &        0&   -0.0002&   -0.0015&    0.0156&    0.0188&   -0.0004\\
         0 &    \cdots&        0 &        0&         0&         0&   -0.0191&    0.1670&    0.1907\\
\end{array} \right|,
 \end{align}
 where constants $T_i(i=1,2,\cdots, 19)$ are given in the Appendix A, and the Matlab code for deriving (\ref{aaa}) is
listed in Appendix B.

It is clear that the last matrix is strictly diagonally dominant, thus its determinant is not equal to zero, which implies that the coefficient matrix of the  CCD system is non-singular. Thus, the CCD linear system (3.16) is uniquely  solvable. This completes the proof of the theorem.
\qed

\section{Numerical experiments}\label{sec5}
    In this section, four numerical experiments are provided to verify the high accuracy of the proposed CCD-TVD method.
All numerical experiments were carried out under MATLAB 2016a on a desktop with 4.00GHz Intel i7-4790K and 16GB RAM.

\begin{example}\label{ex1}
  In this example, we consider the 1D coupled Burgers' equation \cite{LiaoWenyuan.{2008}}:
\begin{equation*}
u_t+uu_x=\frac{1}{Re} u_{xx}, \quad x \in \Omega=[0,1],\quad t\in [0, T],
\end{equation*}
with the following initial condition
\[u(x,0)=\sin(\pi x),\quad \quad\quad x\in (0,1),\]
and boundary conditions
\[u(0,t)=u(1,t)=0,\quad\quad\quad t>0.\]
The exact solution is given by
\[u(x,t)=2\pi \frac{1}{Re}\frac{\Sigma^{\infty}_{n=1}a_{n}\exp(-n^{2}\pi ^{2}\frac{1}{Re} t)n \sin(n\pi x)}{a_{0}+\Sigma^{\infty}_{n=1}a_{n}\exp(-n^{2}\pi ^{2}\frac{1}{Re} t)n \cos(n\pi x)}\]
and the Fourier coefficients
\[a_{0}=\int^{1}_{0}\exp{-(2\pi \frac{1}{Re})^{-1}[1-\cos(\pi x)]}dx,\]

\[a_{n}=2\int^{1}_{0}\exp{-(2\pi \frac{1}{Re})^{-1}[1-\cos(\pi x)]}\cos(n\pi x)dx,n\geq 1.\]

\end{example}

In this problem, we set $v=1, \Delta t=0.00001$, $\Delta x=0.0125$ for all numerical methods. Table \ref{tab1} lists our results and results obtained by researchers~\cite{MuratSari.{2009},MG.{2005},MG.{2006}} using the Hopf-Cole(HC), Restrictive Hopf-Cole(RHC), Restrictive Pad$\acute{e}$ Approximation (RPA) and TVD-CFD (TVCF) methods. It is shown that the CCD-TVD method offers better results than other existing schemes.

\begin{table}[!tbp]
\tabcolsep=10pt
\caption{comparison of results for example 1 at different times for $\frac{1}{Re}=1$, $\Delta t=0.00001, h=0.0125.$}\label{tab1} \centering
  \begin{tabular}{|c|c|c|c|c|c|c|c|}
\hline
  $x$    &  $t$  &  HC\cite{MG.{2005}}  &  RHC\cite{MG.{2005}}  &  RPA\cite{MG.{2006}}  &  TVCF\cite{MuratSari.{2009}}  &  CCD-TVD  &  Exact  \\
\hline  0.25&0.40&0.308860&0.317062&0.308776&0.308894&0.308893&0.308893\\
\hline  0.25&0.60&0.240703&0.248472&0.240739&0.240739&0.240739&0.240739\\
\hline  0.25&0.80&0.195693&0.202953&0.195676&0.195676&0.195676&0.195676\\
\hline  0.25&1.00&0.162561&0.169527&0.162513&0.162565&0.162564&0.162564\\
\hline  0.50&0.40&0.569602&0.583408&0.569527&0.569632&0.569632&0.569632\\
\hline  0.50&0.60&0.447123&0.461714&0.447117&0.447206&0.447205&0.447205\\
\hline  0.50&0.80&0.359152&0.373800&0.359161&0.359236&0.359236&0.359236\\
\hline  0.50&1.00&0.291961&0.306184&0.292843&0.291916&0.291916&0.291916\\
\hline  0.75&0.40&0.625460&0.638847&0.625341&0.625438&0.625437&0.625437\\
\hline  0.75&0.60&0.487337&0.506429&0.487089&0.487215&0.487211&0.487211\\
\hline  0.75&0.80&0.374067&0.393565&0.373827&0.373922&0.373923&0.373923\\
\hline  0.75&1.00&0.287525&0.305862&0.287396&0.287474&0.287473&0.287473\\
\hline
\end{tabular}
\end{table}

\begin{example}\label{ex2}
  In this example, we solve the coupled  2D Burgers' equations~\cite{Liao3, ChenBuyun.{2017}}
\begin{equation*}
u_t+uu_x+vu_y=\frac{1}{Re} (u_{xx}+u_{yy}), \quad (x,y)\in \Omega=[0,1]\times[0,1],\quad t\in [0, T],
\end{equation*}
\begin{equation*}
v_t+uv_x+vv_y=\frac{1}{Re}(v_{xx}+v_{yy}), \quad (x,y)\in \Omega=[0,1]\times[0,1],\quad t\in [0, T],
\end{equation*}
with the following exact solutions:
\begin{equation*}
u(x,y,t)=-\frac{4\pi\frac{1}{Re} e^{-5\pi^2\frac{1}{Re} t} \cos(2\pi x)\sin(\pi y)}{2+e^{-5\pi^2\frac{1}{Re} t}\sin(2\pi x)\sin(\pi y)},
\end{equation*}
\begin{equation*}
v(x,y,t)=-\frac{2\pi\frac{1}{Re} e^{-5\pi^2\frac{1}{Re} t} \sin(2\pi x)\cos(\pi y)}{2+e^{-5\pi^2\frac{1}{Re} t}\sin(2\pi x)\sin(\pi y)}.
\end{equation*}
The initial and boundary conditions are taken from the exact solutions. In this example, we set $\frac{1}{Re}=0.1$.
\end{example}

To carry out the convergence study, we fix $\Delta t=\Delta y^2=\Delta x^2\triangleq h^2$, and change the mesh size $h$ from $1/16$ to $1/256$ by a factor $1/2$.
Results in Table \ref{tab2} show that the current CCD-TVD method is third-order accurate for time variable, and sixth-order accurate for space variable.


\begin{table}[!tbp]
\tabcolsep=10pt
\caption{Numerical results of example 2 by a spatial sixth-order TVD-CCD scheme at T=1 for $\varepsilon$=0.1.}\label{tab2} \centering
  \begin{tabular}{|c|c|c|c|c|}
\hline
$  h    $&$\parallel e_u\parallel_{L^{\infty}}$&Rate& $\parallel e_v\parallel_{L^{\infty}}$&Rate \\
\hline  $1/16 $&5.96E$-$04&  $-$ &1.87E$-$05&  $-$  \\
\hline  $1/32 $&1.92E$-$05&4.96&4.38E$-$07&5.42\\
\hline  $1/64 $&3.04E$-$07&5.98&1.28E$-$08&5.01\\
\hline  $1/128$&5.67E$-$09&5.75&4.19E$-$10&4.93\\
\hline  $1/256$&1.26E$-$10&5.49&1.35E$-$11&4.95\\
\hline
\end{tabular}
\end{table}

\begin{example}\label{ex3}
  In this example, we consider the 2D Burgers' equations, with the initial conditions $$u(x,y,0)=x+y, v(x,y,0)=x-y,$$ and the following exact solutions~\cite{ZHSHDM.{2010}}
\begin{equation*}
u(x,y,z,t)=\frac{x+y-2xt}{1-2t^{2}},
\end{equation*}
\begin{equation*}
v(x,y,z,t)=\frac{x-y-2yt}{1-2t^{2}},
\end{equation*}
where the domain is taken to  be $[0,0.5]\times[0, 0.5]$.  And in this example, we solve the problem at $T=0.1$ with $\frac{1}{Re} =0.1$.
\end{example}


\begin{table}[!tbp]
\tabcolsep=10pt
\caption{Numerical results of example 3 by a spatial sixth-order TVD-CCD scheme at T=0.1 for $\varepsilon$=0.1.}\label{tab3} \centering
  \begin{tabular}{|c|c|c|c|c|}
\hline
$  h    $&$\parallel e_u\parallel_{L^{\infty}}$&Rate& $\parallel e_v\parallel_{L^{\infty}}$&Rate \\
\hline  $1/4 $&1.21E$-$08&  $-$ &9.52E$-$09&  $-$  \\
\hline  $1/8 $&2.10E$-$10&5.85&1.63E$-$10&5.87\\
\hline  $1/16$&3.29E$-$12&5.99&2.59E$-$12&5.98\\
\hline  $1/32$&6.11E$-$14&5.75&4.61E$-$14&5.81\\
\hline
\end{tabular}
\end{table}

\begin{example}\label{ex4}
  In this example, we consider the 3D coupled Burgers' equation (\ref{equation1})-(\ref{equation3}) with following exact solutions:
\begin{equation*}
u(x,y,z,t)=\frac{x+y+z}{1+3t^{2}},
\end{equation*}
\begin{equation*}
v(x,y,z,t)=\frac{x+y+z}{1+3t^{2}},
\end{equation*}
\begin{equation*}
w(x,y,z,t)=\frac{x+y+z}{1+3t^{2}},
\end{equation*}
where the domain is taken to  be $[0,1]\times[0, 1]\times[0, 1]$. The initial and boundary conditions are obtained from the above analytic solutions.  And in this example, we set $\frac{1}{Re}$=0.08.
\end{example}


\begin{table}[!tbp]
\tabcolsep=10pt
\caption{Numerical results of example 4 by a spatial sixth-order TVD-CCD scheme at T=1 for $\varepsilon$=0.08.}\label{tab4} \centering
  \begin{tabular}{|c|c|c|c|c|c|c|}
\hline
$  h    $&$\parallel e_u\parallel_{L^{\infty}}$&Rate& $\parallel e_v\parallel_{L^{\infty}}$&Rate &$\parallel e_w\parallel_{L^{\infty}}$&Rate\\
\hline  $1/4 $&4.84E$-$05&  $-$ &4.84E$-$05&  $-$ &4.84E$-$05&  $-$   \\
\hline  $1/8 $&6.65E$-$07&6.21&6.65E$-$07&6.21&6.65E$-$07&6.21\\
\hline  $1/16$&9.99E$-$09&6.04&9.99E$-$09&6.04&9.99E$-$09&6.04\\
\hline  $1/32$&1.57E$-$10&6.00&1.57E$-$10&6.00&1.57E$-$10&6.00\\
\hline  $1/64$&2.46E$-$12&5.99&2.46E$-$12&5.99&2.46E$-$12&5.99\\
\hline
\end{tabular}
\end{table}

Table \ref{tab3} and Table \ref{tab4} present the results obtained by the CCD-TVD method for Example \ref{ex3} and Example \ref{ex4}, respectively.
From these two tables, it is also shown that the spatial accuracy is sixth-order, and the temporal accuracy is third-order.

\section{Conclusions}\label{sec6}
An efficient spatial sixth-order CCD-TVD method has been developed to solve the multi-dimensional Burger's equation in this paper. This method is verified to be sixth-order accurate in space variable and third-order accurate in time variable. Moreover, the unique solvability of the CCD system is obtained under non-periodic boundary conditions. Existing theoretical works on the CCD method are focused on solving the differential equations subject to periodic boundary conditions~\cite{Wang2018,Gao2}. The convergence analysis for the CCD scheme will be our future objective.

\section*{Acknowledgement}
Kejia Pan was supported by the Natural Science
Foundation of China (No. 41474103), the Excellent Youth Foundation of Hunan Province of China (No. 2018JJ1042)
and the Innovation-Driven Project of Central South University (No. 2018CX042).

\section*{Appendix A.}

       $$T_1= -\frac{136835}{8209824}+\frac{421733\sqrt{7}}{36944208} ,$$
       $$T_2= -\frac{416144963942525\sqrt{7}}{864691128455135232}  ,$$
       $$T_3= -\frac{46840306656146665409\sqrt{7}}{41595480345574524321792}+\frac{6115}{2052456} ,$$
       $$T_4= -\frac{46840306656146665409\sqrt{7}}{41595480345574524321792}+\frac{6115}{2052456} ,$$
       $$T_5= \frac{21881309630676858473952943462656048141172736\sqrt{7}}{4667640113605791995493311956355581953955543731}-\frac{2893014312251833953326629616913521171234816}{518626679289532443943701328483953550439504859} ,$$
       $$T_6=  -\frac{267591658885243284604171740430098010027602763}{33192107474530076412396885022973027228128310976}+\frac{587650275369111747907115169185577623944244693\sqrt{7}}{37341120908846335963946495650844655631644349848} ,$$
       $$T_7=  -\frac{1799614987336256538927773374693436599301849447\sqrt{7}}{298728967270770687711571965206757245053154798784}-\frac{442381615984325718712039486584685244485625}{691502239052709925258268437978604733919339812} ,$$
       $$T_8= \frac{60812707732120381749986704242152551792357469\sqrt{7}}{18670560454423167981973247825422327815822174924}+\frac{900723013413257827633193252372996530470617}{922002985403613233677691250638139645225786416} ,$$
       $$T_9= \frac{5235921989442888980950159\sqrt{7}}{183849407004430256490676224}+\frac{41853249163690425155625}{40855423778762279220150272} ,$$
       $$T_{10}= \frac{640001231916311360619949\sqrt{7}}{551548221013290769472028672}+\frac{33366161622715196997673}{40855423778762279220150272} ,$$
       $$T_{11}= \frac{69251\sqrt{7}}{2574720}+\frac{133}{85824} ,$$
       $$T_{12}= \frac{37\sqrt{7}}{3456} ,$$
       $$T_{13}= -\frac{\sqrt{7}}{17280} ,$$
       $$T_{14}= \frac{34711\sqrt{7}}{8582400}+\frac{7073}{1430400} ,$$
       $$T_{15}= -\frac{3197\sqrt{7}}{1029888}+\frac{2315}{85824}  ,$$
       $$T_{16}= -\frac{4733\sqrt{7}}{2574720}+\frac{479}{107280}  ,$$
       $$T_{17}= -\frac{547\sqrt{7}}{80460}-\frac{29}{26820}  ,$$
       $$T_{18}= \frac{2711\sqrt{7}}{16092}-\frac{1495}{5364} ,$$
       $$T_{19}= \frac{547\sqrt{7}}{80460}+\frac{29}{26820} .$$

\section*{Appendix B. Matlab code for deriving (\ref{aaa})}

\lstset{language=Matlab}
\lstset{breaklines}
\lstset{extendedchars=false}
\begin{lstlisting}
clear;  clc
n = 10; a=6*sqrt(7)/7;  b=3*sqrt(7);
a1=zeros(n,n);a2=zeros(n,n);a3=zeros(n,n);a4=zeros(n,n);
for i=2:n-1
    a1(i,i-1)=7/16; a2(i,i-1)=1/16;   a3(i,i-1)=-9/8; a4(i,i-1)=-1/8;
    a1(i,i)=1;  a2(i,i)=0;  a3(i,i)=0;  a4(i,i)=1;
    a1(i,i+1)=7/16; a2(i,i+1)=-1/16;  a3(i,i+1)=9/8; a4(i,i+1)=-1/8;
end
a1(1,1)=14;a2(1,1)=2;a3(1,1)=1;a4(1,1)=0;
a1(1,2)=16;a2(1,2)=-4;a3(1,2)=2;a4(1,2)=-1;
a1(n,n)=14;a2(n,n)=-2;a3(n,n)=1;a4(n,n)=0;
a1(n,n-1)=16;a2(n,n-1)=4;a3(n,n-1)=2;a4(n,n-1)=1;
a1=sym(a1); a2=sym(a2); a3=sym(a3); a4=sym(a4);

a5=[a1+b*a2            ,a2;
    a3+a*a1+b*(a4+a*a2),a4+a*a2];
a5(n+2:2*n-1,:)=a5(n+2:2*n-1,:)/a5(n+2,2);
a5(n+1,:)=a5(n+1,:)-a5(1,:)*a5(n+1,2)/a5(1,2);
a5(n+n,:)=a5(n+n,:)-a5(n,:)*a5(n+n,n-1)/a5(n,n-1);
a5(n+1,:)=a5(n+1,:)/a5(n+1,1);
a5(n+n,:)=a5(n+n,:)/a5(n+n,n);

a6=a5(1:n,1:n)*a5(n+1:n+n,n+1:n+n)-a5(1:n,n+1:n+n)*a5(n+1:n+n,1:n);
mm11=simplify(a6(1,2)/(a6(2,2)+a6(3,2)*((a6(2,3)*a6(1,2)-a6(2,2)*a6(1,3))/(a6(3,2)*a6(1,3)-a6(3,3)*a6(1,2)))));
mm12=simplify(((a6(2,3)*a6(1,2)-a6(2,2)*a6(1,3))/(a6(3,2)*a6(1,3)-a6(3,3)*a6(1,2)))*mm11);
mm13=(18345*7^(1/2))/1197266 + 405917/513114;
a6(1,:)=a6(1,:)-(mm11*a6(2,:)+mm12*a6(3,:)+mm13*a6(4,:));
m31 = simplify(a6(1,3)/a6(1,1));
a6(:,3)=a6(:,3)-a6(:,1)*m31;
m51 = simplify(a6(1,5)/a6(1,1));
a6(:,5)=a6(:,5)-a6(:,1)*m51;
m23 = simplify(a6(2,3)/a6(3,3));
a6(2,:)=a6(2,:)-a6(3,:)*m23;
a6(:,n-1)=a6(:,n-1)-a6(:,n)*1.5;
mm = simplify(((a6(n,n-2)/a6(n,n))+1/10));
a6(:,n-2)=a6(:,n-2)-a6(:,n)*mm;
a6(n-1,:)=a6(n-1,:)-(a6(n,:)+a6(n-2,:))*0.1;
a6 = simple(a6)

err=double(2*diag(a6) - sum(abs(a6),2))
\end{lstlisting}


\section*{References}
\begin{thebibliography}{45}
\bibitem{Arpad}
{\sc A. Taka$\check{\rm c}$i},{ Mathematical and simulation models of traffic flow}.  PAMM 5 (2005), pp. 633--634.


\bibitem{Yang}
{\sc L. Yang and  X. Pu},  { Derivation of the Burgers' equation from the gas dynamics}. Commun. Math. Sci. 14 (2016), pp. 671--682.


\bibitem{Khesin}
{B. Khesin and G. Misiolek}, { Shock waves for the Burgers equation and curvatures of diffeomorphism groups}.  P. Steklov. I. Math+  259 (2007), pp. 73--81.


\bibitem{Davidson}	
{\sc G. A. Davidson},   { A Burgers' equation approach to finite amplitude acoustics in aerosol media}.  J. Sound. Vib. 38 (1975), pp. 475--495.

\bibitem{Su}	
{\sc N. H. Su, P. C. Watt, K. W. Vincent, M. E. Close, R. Z. Mau},  { Analysis of turbulent flow patterns of soil water under filed conditions using Burgers' equation and porous suction-cup samplers}.  Aust. J. Soil. Res. 42 (2004), pp:9--16.

\bibitem{Zabusky}	
{\sc N. J. Zabusky and  M. D. Kruskal}, { Interaction of solitons in a collisionless plasma and the recurrence of initial states}. { Phys. Rev. Lett.} 15 (1965), pp. 240--243.

\bibitem{Zhao}	
{\sc P. Zhao and M. Qin},  { Multisymplectic geometry and multisysmplectic preissmann scheme for the KdV equation}.  J Phys A-Math Theor. 33 (2000),  pp. 3613--3626.

\bibitem{Esipov}
{\sc S. E. Esipov},  { Coupled Burgers' equations: a model of poly-dispersive sedimentation}. Phys. Rev. E  52 (1995), pp. 3711--3718.

\bibitem{Logan}
{\sc J. D. Logan},  { An introduction to nonlinear partial differential equations} Wily-Interscience, New York, 1994.

\bibitem{Shandarin}
{\sc S. F. Shandarin},  { Three dimensional Burgers' equation as a model for the Large-scale structure Formation in the Universe}. IMA.   85 (1997), pp. 401--413.


\bibitem{Bateman}
{\sc H. Bateman}, { Some recent researches on the motion of fluids}.  Mon. Weather Rev.  43 (1914), pp. 163--170.


\bibitem{Hopf}
{\sc E. Hopf}, { The partial differential equation $u_t+uu_x=\mu u_{xx}$}.   Common. Pur. Appl. Math.  3 (1950), pp. 201--230.

\bibitem{Cole}
{\sc J. D. Cole},  { On a quasilinear parabolic equation occurring in aerodynamics}.  Quart.  Appl.  Math.  9 (1951), pp. 225--236.





\bibitem{Fletcher}
{\sc C. A. Fletcher},  { Generating exact solutions of the two-dimensional Burgers' equations}. Int. J. Numer. Meth. Fl. 3 (1983), pp. 213--216.


\bibitem{GaoZou}
{\sc Q. Gao and M. Zou},  { An analytical solution for two and three dimensional nonlinear Burgers' equation}.  Appl.  Math. Model. 45 (2017), pp. 255-270.







\bibitem{Kutluay2013}
{\sc S. Kutluay and Y. Ucar},  { Numerical solutions of the coupled Burgers' equation by the Galerkin quadratic B-spline finite element method}.  Math. Method. Appl. Sci.  36 (2013), pp. 2403--2415.

\bibitem{Finn}
{\sc E. Var$\ddot{\rm o}$glu and W.D.L. Finn},  { Space-time finite elements incorporating characteristics for the Burgers' equation}. Int. J. Numer. Meth. Eng.  16 (1980), pp. 171--184.


\bibitem{Caldwell}
{\sc J. Caldwell,  P. Wanless and A.E. Cook},  { A finite element approach to Burgers' equation}. Appl.  Math. Model.  5 (1981), pp. 189--193.

\bibitem{Dehghan2014}
{\sc M. Dehghan and  B.N. Saray,  M. Lakestani}  { Mixed finite difference and Galerkin methods for solving Burgers equations using interpolating scaling functions}. Math. Method. Appl. Sci.  37 (2014), pp. 894--912.

\bibitem{Saglam}
{\sc A. R. Bahadir and  M. Sa$\check{\rm g}$lam }, { A mixed finite difference and boundary element approach to one-dimensional Burgers' equation}. Appl. Math. Comput.  160 (2005), pp.  663--673.

\bibitem{Sari}
{\sc M. Sari and G. G$\ddot{\rm u}$rarslan},  { A sixth-order compact finite difference scheme to the numerical solutions of Burgers¡¯ equation}. Appl. Math. Comput. 208 (2009), pp. 475--483.

\bibitem{Kutluay1999}
{\sc S. Kutluay and A. R. Bahadir},  { Numerical solution of one-dimensional Burgers' equation: explicit and exact-explicit finite difference methods}. J.  Comput.  Appl. Math. 103 (1999), pp. 251--261.


\bibitem{Liao1}
{\sc W. Liao },  { An implicit fourth-order compact finite difference scheme for one-dimensional Burgers' equation}.  Appl. Math. Comput. 206 (2008), pp. 755--764.

\bibitem{Liao2}
{\sc W. Liao and J. Zhu}, { Efficient and accurate finite difference schemes for solving one-dimensional Burgers' equation}. Int. J. Comput. Math. 88 (2011), pp. 2575--2590.


\bibitem{Huang}
{\sc P. Huang and A. Abduwali}, { The Modified Local Crank-Nicolson method for one- and two-dimensional Burgers' equations}. Comput. Math. Appl. 59 (2010), pp. 2452--2463.




\bibitem{Bahadir}
{\sc A. R. Bahadir},  { A fully implicit finite-difference scheme for two-dimensional Burgers'  equations}.  Appl. Math. Comput. 137 (2003), pp. 131--137.


\bibitem{Liao3}
{\sc W. Liao}, { A fourth-order finite-difference method for solving the system of two-dimensional Burgers' equations}. Int. J. Numer. Meth. Fl.  64 (2010), pp. 565--590.




\bibitem{Campos}
{\sc M. D. Campos and E. C. Rom$\tilde{\rm a}$o}, { A high-order finite-difference scheme with a linearization technique for solving of three-dimensional Burgers' equation}.  CMES-Comp. Model. Eng.  103 (2014), pp. 139--154.


\bibitem{ChenBuyun.{2017}}
Buyun Chen,Dongdong He,Keijia Pan,A linearized high-order Combined Compact Difference Scheme for Multidimensional Coupled Burgers' equation.Numer.Math.Theor.Meth.Appl.  2(2017) 299-320.

\bibitem{PChu}
{\sc P. Chu and C. Fan}, {  A Three-Point Combined Compact Difference Scheme}. J. Comput. Phys.  140 (1998), pp. 370--399.

\bibitem{Mahesh1998}
{\sc K. Mahesh}, { A family of high order finite difference schemes with good spectral resolution}. J. Comput. Phys. 145 (1998), pp. 332--358.


\bibitem{HSun}
{\sc H. Sun and L. Li}, { A CCD-ADI method for unsteady convection-diffusion equations}. Comput. Phys. Commun. 185 (2014), pp. 790--797.

\bibitem{Nihei2003}
{\sc T. Nihei and  K. Ishii},  { A fast solver of the shallow water equations on a sphere using a combined compact difference scheme}. J. Comput. Phys. 187 (2003), pp. 639--659.

\bibitem{Sun2}
{\sc S. Lee, J. Liu and H. Sun}, {  Combined compact difference scheme for linear second-order
partial differential equations with mixed derivative}. J.  Comput.  Appl. Math. 264 (2014), pp. 23--37.

\bibitem{DDH}
{\sc D. He},  { An unconditionally stable spatial sixth-order CCD-ADI method for the two-dimensional linear hyperbolic equation}.  Numer. Algorithms 72 (2016), pp. 1103--1117.


\bibitem{Gao2}
{\sc G. Gao and H. Sun},  { Three-point combined compact difference schemes for time-fractional advection-diffusion equations with smooth solutions}. J. Comput. Phys. 298  (2015), pp. 520--538.


\bibitem{Li}
{\sc L. Li, H. Sun, S. Tam},  { A spatial sixth-order alternating direction implicit method for two-dimensional cubic nonlinear Schr$\ddot{\rm o}$dinger equations}. Comput. Phys. Commun. 187 (2015), pp. 38--48.

\bibitem{He2017}
{\sc D. He and K. Pan}, { An unconditionally stable linearized CCD-ADI method for generalized nonlinear Schr$\ddot{\rm o}$dinger equations with variable coefficients in two and three dimensions}. Comput. Math. Appl. 73 (2017), pp. 2360--2374.


\bibitem{He20172}
{{\sc D. He and K. Pan}, { A fifth-order combined compact difference scheme
for the Stokes flow on polar geometries}.  E. Asian J. Appl. Math. 7 (2018), pp. 714-727.}

\bibitem{TVD.{1998}}
{\sc S. Gottlieb and C. Shu},{ Total variation diminishing Runger-Kutta schemes}. Math. Comput. 221 (1998), pp. 73--85.

\bibitem{LiaoWenyuan.{2008}}
{\sc W. Liao},{An implicit fourth-order compact finite difference scheme for one-dimensional Burgers' equation} . Appl. Math. Comput. 206 (2008), pp. 755--764.

\bibitem{MuratSari.{2009}}
{\sc M. Sari and G. Gurarslan} ,{A sixth-order compact finite difference scheme to the numerical solutions of Burgers' equation}.Appl. Math. Comput. 208 (2009), pp. 475--483.


\bibitem{MG.{2005}}
{\sc M.G$\ddot{u}$lsu and T.$\ddot{O}$zis} , {Numerical solution of Burgers' equation with restrictive Taylor approximations} . Appl. Math. Comput. 171(2005), pp. 1192--1200.

\bibitem{MG.{2006}}
{\sc M.G$\ddot{u}$lsu} , {A finite difference approach for solution of Burgers' equation} . Appl. Math. Comput. 175 (2006), pp. 1245--1255.

\bibitem{ZHSHDM.{2010}}
{\sc H. Zhu,H. Shu and M. Ding} , {Numerical solutions of two-dimensional Burgers' equations by discrete Adomian decomposition method} . Comput. Math. Appl. 60 (2010), pp. 840--848.

\bibitem{Wang2018}
{\sc Q.H. Wang, K.J. Pan and H.L. Hu}, {Unique solvability of the CCD scheme for convection¨Cdiffusion equations with variable convection coefficients}.
Adv. Diff. Equ. 2018 (2018), pp. 163.


\end{thebibliography}
\end{document}